\documentclass[12pt,oneside,reqno]{amsart}

\usepackage{amsmath,amssymb,amsthm,textcomp}
\usepackage{amsfonts,graphicx}
\usepackage[mathscr]{eucal}
\usepackage{color}
\usepackage{csquotes}
\usepackage[backend=bibtex,%
firstinits=true,%
doi=false,%
isbn=true,%
url=false,%
maxnames=99]{biblatex}%

\AtEveryBibitem{\clearfield{issn}}
\AtEveryCitekey{\clearfield{issn}}
\numberwithin{equation}{section}
\DeclareNameAlias{sortname}{last-first}
\theoremstyle{definition}
\usepackage{mathtools}
\numberwithin{equation}{section}

\newcommand{\ncom}{\newcommand}

\ncom{\beq}{\begin{equation}}
\ncom{\eeq}{\end{equation}}
\ncom{\bea}{\begin{eqnarray*}}
\ncom{\eea}{\end{eqnarray*}}
\ncom{\beqa}{\begin{eqnarray}}
\ncom{\eeqa}{\end{eqnarray}}
\ncom{\nno}{\nonumber}
\ncom{\non}{\nonumber}
\ncom{\ds}{\displaystyle}
\ncom{\half}{\frac{1}{2}}
\ncom{\mbx}{\makebox{.25cm}}
\ncom{\hs}{\mbox{\hspace{.25cm}}}
\ncom{\rar}{\rightarrow}
\ncom{\Rar}{\Rightarrow}
\ncom{\noin}{\noindent}
\ncom{\bc}{\begin{center}}
\ncom{\ec}{\end{center}}
\ncom{\sz}{\scriptsize}
\ncom{\rf}{\ref}
\ncom{\s}{\sqrt{2}}
\ncom{\sgm}{\sigma}
\ncom{\Sgm}{\Sigma}
\ncom{\psgm}{\sigma^{\prime}}
\ncom{\dt}{\delta}
\ncom{\Dt}{\Delta}
\ncom{\lmd}{\lambda}
\ncom{\Lmd}{\Lambda}
\ncom{\Th}{\Theta}
\ncom{\e}{\eta}
\ncom{\eps}{\epsilon}
\ncom{\pcc}{\stackrel{P}{>}}
\ncom{\lp}{\stackrel{L_{p}}{>}}
\ncom{\dist}{{\rm\,dist}}
\ncom{\sspan}{{\rm\,span}}
\ncom{\re}{{\rm Re\,}}
\ncom{\im}{{\rm Im\,}}
\ncom{\sgn}{{\rm sgn\,}}
\ncom{\ba}{\begin{array}}
\ncom{\ea}{\end{array}}
\ncom{\hone}{\mbox{\hspace{1em}}}
\ncom{\htwo}{\mbox{\hspace{2em}}}
\ncom{\hthree}{\mbox{\hspace{3em}}}
\ncom{\hfour}{\mbox{\hspace{4em}}}
\ncom{\vone}{\vskip 2ex}
\ncom{\vtwo}{\vskip 4ex}
\ncom{\vonee}{\vskip 1.5ex}
\ncom{\vthree}{\vskip 6ex}
\ncom{\vfour}{\vspace*{8ex}}
\ncom{\norm}{\|\;\;\|}
\ncom{\integ}[4]{\int_{#1}^{#2}\,{#3}\,d{#4}}
\ncom{\vspan}[1]{{{\rm\,span}\{ #1 \}}}
\ncom{\dm}[1]{ {\displaystyle{#1} } }
\ncom{\ri}[1]{{#1} \index{#1}}

\newtheorem{remark}{\bf Remark}[section]

\newtheoremstyle
    {remarkstyle}
    {}
    {11pt}
    {}
    {}
    {\bfseries}
    {:}
    {     }
    {\thmname{#1} \thmnumber{#2} }

\theoremstyle{remarkstyle}

\def\eps{\varepsilon}

\begin{document}
\title{On Elephant Random Walk with Random Memory}
\author[Manisha Dhillon]{Manisha Dhillon}
\address{Manisha Dhillon, Department of Mathematics, Indian Institute of Technology Bhilai, Durg 491002, India.}
\email{manishadh@iitbhilai.ac.in}
\author[Kuldeep Kumar Kataria]{Kuldeep Kumar Kataria}
\address{Kuldeep Kumar Kataria, Department of Mathematics, Indian Institute of Technology Bhilai, Durg 491002, India.}
\email{kuldeepk@iitbhilai.ac.in}
\date{\today}
\begin{abstract}
 In this paper, we introduce the elephant random walk (ERW) with memory consisting of randomly selected steps from its history. It is a time-changed variant of the standard  elephant random walk with memory consisting of its full history. At each time point, the time changing component is the composition of two uniformly distributed independent random variables with support over all the past steps. Several conditional distributional properties including the conditional mean increments and conditional displacement of ERW with random memory are obtained. Using these conditional results, we derive the recursive and explicit expressions for the mean increments and mean displacement of the walk.
\end{abstract}

\maketitle
\section{Introduction}
Random walks have wide range of applications in the fields such as biosciences (see Velleman (2014)), order statistics (see Blondel {\it et al.} (2020)), econometrics (see Gourieroux and Jasiak (2022)), {\it etc.} The simple symmetric random walk is one of the most extensively studied models among various types of random walks. In simple symmetric random walk, the steps are of unit size that are chosen with equal probability and are independent of each other. That is, the walker has no memory. This random walk exhibits the diffusive behaviour. However, anomalous diffusion is observed in many physical and biological systems in which the theoretical models incorporate memory effects from the previous steps. The elephant random walk (ERW) is one among the random walks with anomalous diffusion. It was introduced by Schütz and Trimper (2004) to study the effect of memory on random walk. In case of standard ERW, the walker takes into account the complete memory, that is, the next step depends on all the past steps taken so far. Its name is inspired by the fact that elephants have long memory. In the past two decades, the ERW has gained attention of several researchers, for example, Bercu (2018), Bertoin (2022),   {\it etc}.

First, we briefly describe the standard ERW. It is a one-dimensional discrete-time random walk on $\mathbb{Z}$ such that starting from the origin it performs one step of unit size at each time point. At the first time point, the walker takes a unit step towards right with probability $q\in[0,1]$ or a unit step towards left with probability $1-q$. For each subsequent steps, the walker chooses one of its previous steps and either repeats it with probability $p\in[0,1]$ or moves in the opposite direction with probability $1-p$. Thus, the first step $X_1$ is distributed as
\begin{equation*}
	X_{1}=\begin{cases}
		+1\ \ \text{with probability $q$},\\
		-1\ \ \text{with probability $1-q$}\\
	\end{cases}
\end{equation*}
and after $n$ steps, that is, at position $S_n=\sum_{i=1}^{n}X_i$, $n\ge1$ we have
\begin{equation*}
	X_{n+1}=\begin{cases}
		+X_{K}\ \ \text{with probability $p$},\\
		-X_{K}\ \ \text{with probability $1-p$},
	\end{cases}
\end{equation*}
where $K\sim$ unif$\{1,2,\dots,n\}$. It is important to note that the ERW exhibits a time-inhomogeneous Markovian structure.

Several extensions of the standard ERW have been studied in literature, for example, Baur and Bertoin (2016), Gut and Stadtm\"uller (2021), Fan and Shao (2024), \textit{etc.}
Some refinements on the asymptotic behaviour of one-dimensional ERW using martingale approach are obtained by Bercu (2018). A few results on the law of large numbers, central limit theorem and strong invariance principle for the standard ERW are derived by Coletti \textit{et al.} (2017a), (2017b).
Gut and Stadtm\"uller (2021) study the ERW with restricted memory, for instance, the walker remembers only some distant past, only a recent past, or a mixture of both.  For the ERW with random step sizes, we refer the reader to Fan and Shao (2024), Dedecker \textit{et al.} (2023), Roy \textit{et al.} (2025), and the references therein. 

In this paper, we introduce and study the ERW with random memory. At each time point, the random memory set is determined by the rolls of an unbiased die whose outcomes are independent of the walk. It is a time-changed variant of the standard ERW with memory consisting of its full history. At each time point, the time changing component is the composition of two uniformly distributed independent random variables with support over all the past steps. The conditional distribution of the increments of ERW with random memory is determined. From these conditional results, we establish recursive relations for the mean increments and the mean displacement of the walker. Also, their explicit expressions are obtained.
\section{ERW with random memory}
We consider a one-dimensional elephant random walk $\{S_n\}_{n\ge0}$ on integers with a random memory set. 
The walk starts at the origin, that is, $S_0=0$. At $n=1$, the walker moves a unit step towards right with probability $q\in[0,1]$ and a unit step towards left with probability $1-q$. So, the walker's first step $X_1$ is Rademacher $\mathcal{R}(q)$ distributed. At each time step $n\ge1$, a $n$-faced unbiased die is rolled and let $Y(n)\in\{1,2,\dots,n\}$ be the outcome of the roll. Let $K$ be a discrete uniform random variable on $\{1,2,\dots, Y(n)\}$, that is, $K\sim$ unif$\{1,2,\dots,Y(n)\}$. Then, the walker takes the next step according to the following law:
\begin{equation*}
X_{n+1}=\begin{cases}
	+X_{K}\ \ \text{with probability $p$},\\
	-X_{K}\ \ \text{with probability $1-p$},
\end{cases}
\end{equation*}
where $0\le p\le 1$. 
Let $S_{n+1}$ be the position of the walker at time $n+1$.  Then, 
\begin{equation}\label{STEVEQ}
S_{n+1}=S_n+X_{n+1}, \ n\ge 0.
\end{equation} 

For $p=1/2$ and $q=1/2$, the ERW $\{S_n\}_{n\ge0}$ reduces to the simple symmetric random walk.

Note that for any $n\ge1$, the $(n+1)$-th step of walker can be written as $X_{n+1}=\alpha_nX_{\beta(Y(n))}$. Here, $\alpha_n$, $\beta(n)$ and $Y(n)$ are independent discrete random variables such that $\alpha_n$ has a Rademacher $\mathcal{R}(p)$ distribution, $\beta(n)\sim$ unif$\{1,2,\dots, n\}$ and $Y(n)\sim$ unif$\{1,2,\dots,n\}$. Also, $\alpha_n$, $\beta(n)$ and $Y(n)$ are independent of $X_1, X_2,\dots, X_n$.

Let $\mathcal{M}_n=\{1,2,\dots, Y(n)\}\subset \{1,2,\dots,n\}$ be the random memory set of the walker. That is, $\mathcal{M}_n$ contains those steps from the first $n$ steps on which the walker decides the next step. 

For $n\ge1$, let $\mathcal{G}_n=\sigma\{X_1,X_2,\dots,X_n\}$ be the $\sigma$-algebra generated by complete past up to step $n$ and $\mathcal{F}_n=\sigma\{X_k,\, k\in\mathcal{M}_n\}$ be the $\sigma$-algebra generated by the random memory set. 
Then, for $x\in\{-1,1\}$, we have
\begin{align}
\mathbb{P}\{X_{n+1}=x|\mathcal{F}_n\}&=\mathbb{P}\{\alpha_nX_{\beta(Y(n))}=x|\mathcal{F}_n\}\nonumber\\
&=\sum_{r=1}^{n}\mathbb{P}\{Y(n)=r\}\mathbb{P}\{\alpha_nX_{\beta(r)}=x|\mathcal{F}_n, Y(n)=r\} \nonumber\\
&=\sum_{r=1}^{n}\frac{1}{n}\sum_{k=1}^{r}\frac{1}{r}\mathbb{P}\{\alpha_nX_{k}=x|\mathcal{F}_n, Y(n)=r\},\ n\ge1, \label{genpm}
\end{align}
where the penultimate step follows on using the independence of $\beta(n)$ and $Y(n)$.

From \eqref{genpm}, for $x=1$, we get
\begin{align}
\mathbb{P}\{X_{n+1}=1|\mathcal{F}_n\}&=\sum_{r=1}^{n}\frac{1}{n}\sum_{k=1}^{r}\frac{1}{r}\mathbb{P}\{\alpha_nX_{k}=1|\mathcal{F}_n, Y(n)=r\}\nonumber\\
&=\sum_{r=1}^{n}\sum_{k=1}^{r}\frac{1}{nr}\Big(\mathbb{P}\{\alpha_n=1, X_{k}=1|\mathcal{F}_n, Y(n)=r\}\nonumber\\
&\hspace{4cm}+\mathbb{P}\{\alpha_n=-1, X_{k}=-1|\mathcal{F}_n, Y(n)=r\}\Big)\nonumber\\
&=\sum_{r=1}^{n}\sum_{k=1}^{r}\frac{1}{nr}\Big(\mathbb{P}\{\alpha_n=1, X_{k}=1|\sigma(X_1,X_2,\dots,X_r)\}\nonumber\\
&\hspace{3.2cm}+\mathbb{P}\{\alpha_n=-1, X_{k}=-1|\sigma(X_1,X_2,\dots,X_r)\}\Big).\label{pmx1}
\end{align}
For any $G\in \sigma(X_1,X_2,\dots,X_r)$, we have
\begin{align*}
\int_{G}\mathbb{P}\{\alpha_n=1, X_{k}=1|\sigma(X_1,X_2,\dots,X_r)\}\mathrm{d}\mathbb{P}&=\mathbb{P}(\{\alpha_n=1\}\cap  \{X_{k}=1\}\cap G)\\
&=\int_{\{X_{k}=1\}\cap G}\mathbb{P}\{\alpha_n=1\}\mathrm{d}\mathbb{P}\\
&=p\int_{G}I_{\{X_k=1\}}\mathrm{d}\mathbb{P},
\end{align*}
where $I_{A}$ denotes the indicator function on set $A$.
As $\{X_k=1\}\in\mathcal{F}_k$ $\subset$ $\sigma(X_1,X_2,\dots, X_r)$, 
\begin{equation}\label{pmx1a}
\mathbb{P}\{\alpha_n=1, X_{k}=1|\sigma(X_1,X_2,\dots,X_r)\}=pI_{\{X_k=1\}}
\end{equation}
with probability $1$.

Similarly,
\begin{equation}\label{pmx1b}
\mathbb{P}\{\alpha_n=-1, X_{k}=-1|\sigma(X_1,X_2,\dots,X_r)\}=(1-p)I_{\{X_k=-1\}}
\end{equation}
with probability $1$.
By substituting \eqref{pmx1a} and \eqref{pmx1b} in \eqref{pmx1}, we get
\begin{align}\label{pmfgrh}
\mathbb{P}\{X_{n+1}=1|\mathcal{F}_n\}&=\sum_{r=1}^{n}\sum_{k=1}^{r}\frac{1}{nr}\Big(pI_{\{X_k=1\}}+(1-p)I_{\{X_k=-1\}}\Big)\nonumber\\
&=\sum_{r=1}^{n}\sum_{k=1}^{r}\frac{1}{2nr}(1+(2p-1)X_k).
\end{align}

Similarly, for $x=-1$, we have
\begin{align}\label{pmfgrhh}
\mathbb{P}\{X_{n+1}=-1|\mathcal{F}_n\}&=\sum_{r=1}^{n}\sum_{k=1}^{r}\frac{1}{nr}\Big(pI_{\{X_k=-1\}}+(1-p)I_{\{X_k=1\}}\Big)\nonumber\\
&=\sum_{r=1}^{n}\sum_{k=1}^{r}\frac{1}{2nr}(1-(2p-1)X_k).
\end{align}
By using \eqref{pmfgrh} and \eqref{pmfgrhh} in \eqref{genpm}, we get
\begin{equation*}
\mathbb{P}\{X_{n+1}=x|\mathcal{F}_n\}=\sum_{r=1}^{n}\sum_{k=1}^{r}\frac{1}{2nr}(1+(2p-1)xX_k),\ n\ge1.
\end{equation*}
 Thus, the conditional mean is given by
\begin{equation}\label{eqcondi}
\mathbb{E}(X_{n+1}|\mathcal{F}_n)=\sum_{x=\pm1}x\mathbb{P}\{X_{n+1}=x|\mathcal{F}_n\}=\frac{(2p-1)}{n}\sum_{r=1}^{n}\sum_{k=1}^{r}\frac{X_k}{r},\ n\ge1.
\end{equation}

The following recursive relation of mean increments is obtained by using \eqref{eqcondi} and the law of iterated expectations:
\begin{equation}\label{eqcondii}
	\mathbb{E}(X_{n+1})=\frac{(2p-1)}{n}\sum_{r=1}^{n}\sum_{k=1}^{r}\frac{\mathbb{E}(X_k)}{r},\ n\ge1.
\end{equation}
Also, from \eqref{STEVEQ} and \eqref{eqcondi}, we have the following conditional mean of ERW with random memory:
\begin{equation*}
	\mathbb{E}(S_{n+1}|\mathcal{F}_n)=\mathbb{E}(S_{n}|\mathcal{F}_n)+\frac{(2p-1)}{n}\sum_{r=1}^{n}\frac{S_r}{r}, \ n\ge 1.
\end{equation*}
Thus, we get the following recursive relation for its mean:
\begin{equation}\label{meanSn}
\mathbb{E}(S_{n+1})=\mathbb{E}(S_n)+\frac{(2p-1)}{n}\sum_{r=1}^{n}\frac{\mathbb{E}(S_r)}{r}
\end{equation}
which can also be obtained from \eqref{STEVEQ} and \eqref{eqcondii}.
\begin{remark}
If we condition on the steps that are not contained in the random memory set then the walker can not choose them for the immediate next step, that is,
\begin{equation*}
\mathbb{E}(X_{n+1}|\mathcal{G}_n)=\mathbb{E}(X_{n+1}|\mathcal{F}_n)=\frac{(2p-1)}{n}\sum_{r=1}^{n}\sum_{k=1}^{r}\frac{X_k}{r}.
\end{equation*}
\end{remark}
\begin{remark}
We note that the equalities in \eqref{pmx1a} and \eqref{pmx1b} hold true with probability $1$. In all the expressions for conditional probabilities and conditional expectations, the phrase `with probability $1$' is dropped to avoid repetition.
\end{remark}

	Next, we obtain a recursive relation for the mean of $X_{m+1}X_{n+1}$, $m\ge1$, $n\ge1$. 
	
	For $m=n$, we have $\mathbb{E}(X_{n+1}^2)=1$. 	
	For $m\ne n$, let $\mathcal{H}_{m,n}=\sigma(\mathcal{F}_m\cup\mathcal{F}_n)$. Then, for $x=\pm1$, we have 
	\begin{align}\label{pmxixj}
			\mathbb{P}\{X_{m+1}X_{n+1}=x|\mathcal{H}_{m,n}\}&=\mathbb{P}\{\alpha(m)X_{\beta(Y(m))}\alpha(n)X_{\beta(Y(n))}=x|\mathcal{H}_{m,n}\}\nonumber\\
			&=\sum_{r=1}^{n}\sum_{l=1}^{m}\mathbb{P}\{Y(n)=r\}\mathbb{P}\{Y(m)=l\}\nonumber\\
			&\hspace{1.5cm} \cdot\mathbb{P}\{\alpha(m)\alpha(n)X_{\beta(r)}X_{\beta(l)}=x|\mathcal{H}_{m,n},Y(n)=r,Y(m)=l\}\nonumber\\
			&=\sum_{r=1}^{n}\sum_{l=1}^{m}\sum_{k=1}^{r}\sum_{j=1}^{l}\frac{1}{mnrl}\nonumber\\
			&\hspace{.8cm}\cdot\mathbb{P}\{\alpha(m)\alpha(n)X_{k}X_{j}=x|\mathcal{H}_{m,n},Y(n)=r,Y(m)=l\},
	\end{align}
	where we have used the independence of $\beta(n)$ and $Y(n)$ to get the penultimate step.
	
	Let $E_{i_1,i_2,i_3,i_4}=\{\alpha(m)=i_1,\alpha(n)=i_2,X_{k}=i_3,X_{j}=i_4\}$. On taking $x=1$ in \eqref{pmxixj}, we obtain
	{\small\begin{equation}\label{condiXmn1}
			\mathbb{P}\{X_{m+1}X_{n+1}=1|\mathcal{H}_{m,n}\}=\sum_{r=1}^{n}\sum_{l=1}^{m}\sum_{k=1}^{r}\sum_{j=1}^{l}\frac{1}{mnrl}\sum_{\substack{\prod_{y=1}^{4}i_y =1,\\
			i_y\in\{-1,1\}}} \mathbb{P}\{E_{i_1,i_2,i_3,i_4}|\mathcal{H}_{m,n},Y(n)=r,Y(m)=l\}.
	\end{equation}}
	Let $G\in\sigma(\sigma(X_1,X_2,\dots,X_r)\cup\sigma(X_1,X_2,\dots,X_l))$. Then,
	\begin{align}\label{RDth}
		\int_{G}\mathbb{P}\{E_{1,1,1,1}|\mathcal{H}_{m,n},Y(n)=r,Y(m)=l\}\mathrm{d}\mathbb{P}&=\mathbb{P}(E_{1,1,1,1}\cap G)\nonumber\\
		&=\displaystyle\int_{\{X_{k}=1,X_{j}=1\}\cap G}\mathbb{P}\{\alpha(m)=1, \alpha(n)=1\}\mathrm{d}\mathbb{P}\nonumber\\
		&=p^2\int_{G}I_{\{X_{k}=1,X_{j}=1\}}\mathrm{d}\mathbb{P},
	\end{align}
	where the last step follows from the independence of $\alpha(m)$ and $\alpha(n)$.
	
	 Note that $\{X_{k}=1,X_{j}=1\}\in \sigma(\sigma(X_1,X_2,\dots,X_r)\cup\sigma(X_1,X_2,\dots,X_l))$. So, from \eqref{RDth}, we get
	\begin{equation}\label{alcas1}
			\mathbb{P}\{E_{1,1,1,1}|\mathcal{H}_{m,n},Y(n)=r,Y(m)=l\}=p^2I_{\{X_{k}=1,X_{j}=1\}}.
	\end{equation}
	Similarly, we have
	\begin{align}\label{alcas2}
			\left.\begin{aligned}
				\mathbb{P}\{E_{1,1,-1,-1}|\mathcal{H}_{m,n},Y(n)=r,Y(m)=l\}&=p^2I_{\{X_{k}=-1,X_{j}=-1\}},\\
				\mathbb{P}\{E_{1,-1,1,-1}|\mathcal{H}_{m,n},Y(n)=r,Y(m)=l\}&=p(1-p)I_{\{X_{k}=1,X_{j}=-1\}},\\
				\mathbb{P}\{E_{1,-1,-1,1}|\mathcal{H}_{m,n},Y(n)=r,Y(m)=l\}&=p(1-p)I_{\{X_{k}=-1,X_{j}=1\}},\\
				\mathbb{P}\{E_{-1,-1,1,1}|\mathcal{H}_{m,n},Y(n)=r,Y(m)=l\}&=(1-p)^2I_{\{X_{k}=1,X_{j}=1\}},\\
				\mathbb{P}\{E_{-1,1,-1,1}|\mathcal{H}_{m,n},Y(n)=r,Y(m)=l\}&=p(1-p)I_{\{X_{k}=-1,X_{j}=1\}},\\
				\mathbb{P}\{E_{-1,1,1,-1}|\mathcal{H}_{m,n},Y(n)=r,Y(m)=l\}&=p(1-p)I_{\{X_{k}=1,X_{j}=-1\}},\\
				\mathbb{P}\{E_{-1,-1,-1,-1}|\mathcal{H}_{m,n},Y(n)=r,Y(m)=l\}&=(1-p)^2I_{\{X_{k}=-1,X_{j}=-1\}}.
			\end{aligned}
			\right\}
	\end{align}
	By substituting \eqref{alcas1} and \eqref{alcas2} in \eqref{condiXmn1}, we get
	\begin{align}\label{x1xnxm}
			\mathbb{P}&\{X_{m+1}X_{n+1}=1|\mathcal{H}_{m,n}\}\nonumber\\
			&=\sum_{r=1}^{n}\sum_{l=1}^{m}\sum_{k=1}^{r}\sum_{j=1}^{l}\frac{1}{mnrl}\Big(pI_{\{X_{k}=1\}}\Big(\frac{1+(2p-1)X_{j}}{2}\Big) +pI_{\{X_{k}=-1\}}\Big(\frac{1-(2p-1)X_{j}}{2}\Big)\nonumber\\
			&\hspace{2cm} +(1-p)I_{\{X_{k}=-1\}}\Big(\frac{1+(2p-1)X_{j}}{2}\Big)+(1-p)I_{\{X_{k}=1\}}\Big(\frac{1-(2p-1)X_{j}}{2}\Big)\Big)\nonumber\\
			&=\sum_{r=1}^{n}\sum_{l=1}^{m}\sum_{k=1}^{r}\sum_{j=1}^{l}\frac{1}{mnrl}\Big(\Big(\frac{1+(2p-1)X_{j}}{2}\Big)\Big(\frac{1+(2p-1)X_{k}}{2}\Big)\nonumber\\ &\hspace{7cm}+\Big(\frac{1-(2p-1)X_{j}}{2}\Big)\Big(\frac{1-(2p-1)X_{k}}{2}\Big)\Big). 
	\end{align}
	Similarly, for $x=-1$, we have
	{\small\begin{align}\label{x-1xnxm}
			\mathbb{P}\{X_{m+1}X_{n+1}=-1|\mathcal{H}_{m,n}\}&=\sum_{r=1}^{n}\sum_{l=1}^{m}\sum_{k=1}^{r}\sum_{j=1}^{l}\frac{1}{mnrl}\Big(\Big(\frac{1-(2p-1)X_{j}}{2}\Big)\Big(\frac{1+(2p-1)X_{k}}{2}\Big)\nonumber\\
			&\hspace{3cm} +\Big(\frac{1+(2p-1)X_{j}}{2}\Big)\Big(\frac{1-(2p-1)X_{k}}{2}\Big)\Big).
	\end{align}}
	From \eqref{x1xnxm} and \eqref{x-1xnxm}, we get
\begin{equation*}
\mathbb{E}(X_{m+1}X_{n+1}|\mathcal{H}_{m,n})=(2p-1)^2\sum_{r=1}^{n}\sum_{l=1}^{m}\sum_{k=1}^{r}\sum_{j=1}^{l}\frac{1}{mnrl}X_jX_k.
\end{equation*}
Thus, we have the following recursive relation:
\begin{equation}\label{meanXmXnexp}
\mathbb{E}(X_{m+1}X_{n+1})=
	(2p-1)^2\sum_{r=1}^{n}\sum_{l=1}^{m}\sum_{k=1}^{r}\sum_{j=1}^{l}\frac{1}{mnrl}\mathbb{E}(X_jX_k),\ m\ne n, \, m\ge1,\, n\ge1.
\end{equation}

From \eqref{STEVEQ} and \eqref{meanXmXnexp}, we get the second moment of ERW with random memory in the following form:
\begin{equation*}
\mathbb{E}(S_{n+1}^2)=n+1+2(2p-1)^2\sum_{0\le i< j\le n}\sum_{r=1}^{i}\sum_{l=1}^{j}\sum_{k=1}^{r}\sum_{q=1}^{l}\frac{\mathbb{E}(X_kX_q)}{ijrl}, \ n\ge0.
\end{equation*}

\section{Mean of ERW with random memory}
Here, we derive the explicit expressions for the mean increments and mean displacement of the walker performing ERW with random memory. 

First, we obtain the expressions for mean increments.
\subsection{Mean increments}
Let $X_1=1$ and $\alpha=2p-1$. Then, by using \eqref{eqcondii}, we have
\begin{align*}
	\mathbb{E}(X_{n+1})=\frac{\alpha}{n}\sum_{r=1}^{n}\sum_{k=1}^{r}\frac{\mathbb{E}(X_k)}{r}.
\end{align*}
So, $\mathbb{E}(X_1)=1$, $\mathbb{E}(X_2)=\alpha\mathbb{E}(X_1)=\alpha$ and 
\begin{equation*}
	\mathbb{E}(X_3)=\frac{\alpha}{2}\Big(\mathbb{E}(X_1)+\frac{\mathbb{E}(X_1)+\mathbb{E}(X_2)}{2}\Big)=\frac{\alpha}{2}\Big(1+\frac{1}{2}\Big)+\frac{\alpha^2}{2^2}=\sum_{j=1}^{2}\alpha^j\sum_{x_2\in\Phi^j_2}\frac{1}{2^{x_2}},
\end{equation*}
where $\Phi^1_2=\{1, 2\}$ and $\Phi^2_2=\{2\}$. 

For $n=3$ and $n=4$, we have
\begin{align*}
	\mathbb{E}(X_4)&=\frac{\alpha}{3}\Big(\mathbb{E}(X_1)+\frac{\mathbb{E}(X_1)+\mathbb{E}(X_2)}{2}+\frac{\mathbb{E}(X_1)+\mathbb{E}(X_2)+\mathbb{E}(X_3)}{3}\Big)\\
	&=\frac{\alpha}{3}\Big(1+\frac{1}{2}+\frac{1}{3}\Big)+\alpha^2\Big(\frac{1}{2\cdot3}+\frac{1}{3\cdot3}+\frac{1}{2\cdot3^2}+\frac{1}{2^2\cdot3^2}\Big)+\frac{\alpha^3}{2^2\cdot3^2}\\
	&=\sum_{j=1}^{3}\alpha^j\sum_{(x_2,x_3)\in\Phi^j_3}\frac{1}{2^{x_2}3^{x_3}},\\
	\mathbb{E}(X_5)&=\frac{\alpha}{4}\Big(\mathbb{E}(X_1)+\frac{\mathbb{E}(X_1)+\mathbb{E}(X_2)}{2}\\
	&\hspace{1.8cm} +\frac{\mathbb{E}(X_1)+\mathbb{E}(X_2)+\mathbb{E}(X_3)}{3}+\frac{\mathbb{E}(X_1)+\mathbb{E}(X_2)+\mathbb{E}(X_3)+\mathbb{E}(X_4)}{4}\Big)\\
	&=\frac{\alpha}{4}\Big(1+\frac{1}{2}+\frac{1}{3}+\frac{1}{4}\Big)+\alpha^2\Big(\frac{1}{2\cdot4}+\frac{1}{2\cdot3\cdot4}+\frac{1}{2^2\cdot3\cdot4}+\frac{1}{3\cdot4}+\frac{1}{4^2}+\frac{1}{2\cdot4^2}\\
	&\ \ +\frac{1}{2^2\cdot 4^2}+\frac{1}{3\cdot4^2}+\frac{1}{2\cdot3\cdot4^2}+\frac{1}{3^2\cdot4^2}\Big)+\alpha^3\Big(\frac{1}{2^2\cdot3\cdot4}+\frac{1}{2^2\cdot4^2}+\frac{1}{2\cdot3\cdot4^2}\\
	&\ \ +\frac{1}{3^2\cdot4^2}+\frac{1}{2\cdot3^2\cdot4^2}+\frac{1}{2^2\cdot3^2\cdot4^2}\Big)+\frac{\alpha^4}{2^2\cdot3^2\cdot4^2}\\
	&=\sum_{j=1}^{4}\alpha^j\sum_{(x_2,x_3, x_4)\in\Phi^j_4}\frac{1}{2^{x_2}3^{x_3}4^{x_4}},
\end{align*} 
where $\Phi^1_3=\{(0,1),(1,1),(0,2)\}$, $\Phi^2_3=\{(1,1),(0,2),(1,2,),(2,2)\}$, $\Phi^3_3=\{(2,2)\}$, $\Phi^1_4=\{(1,0,1),(0,1,1),(0,0,1),(0,0,2)\}$, $\Phi^2_4=\{(1,0,1),(1,1,1),(2,1,1),(0,1,1),(0,0,2),  (1,0,2)$, $  (2,0,2),(0,1,2),(1,1,2),(0,2,2)\}$, $\Phi^3_4=\{(2,1,1),(2,0,2),(1,1,2),(0,2,2),(1,2,2),(2,2,2)\}$ and $\Phi^4_4=\{(2,2,2)\}$.

Proceeding inductively, we get the mean increments of ERW with random memory in the following form:
\begin{equation}\label{Ex_n+1}
	\mathbb{E}(X_{n+1})=\sum_{j=1}^{n}\alpha^j\sum_{(x_2,x_3,\dots,x_n)\in\Phi^j_n}\frac{1}{2^{x_2}3^{x_3}\dots n^{x_n}}, \ n\ge2.
\end{equation}
Here, $\Phi^1_n=\{(0,0,\dots,0,2), (0,0,\dots,0,1), (0,0,\dots,0,1,1), (0,0,\dots,0,1,0,1),\dots,(0,1,0$, $\dots, 0,1),(1,0,\dots,0,1)\}\subset\mathbb{R}^{n-1}$ and 
\begin{align*}
	\Phi^j_n=\begin{cases}
		(\mathcal{B}^{1}_{n-1}\cup\mathcal{B}^{1}_{n-2}\cup\dots\cup\mathcal{B}^{1}_{2})\cup \Phi^{1}_n\backslash\{(0,0,\dots,0,1)\},\ j=2,\\
		\mathcal{B}^{j-1}_{n-1}\cup\mathcal{B}^{j-1}_{n-2}\cup\dots\cup\mathcal{B}^{j-1}_{j-1},\ j\ge3,
	\end{cases}
\end{align*} 
where for $j\ge 2$ and $1\le k\le n-j+1$, $\mathcal{B}^{j-1}_{n-k}$'s are
{\small\begin{align*}
		\mathcal{B}^{j-1}_{n-1}&=\{(x_2,x_3,\dots,x_{n-1},2):(x_2,x_3,\dots,x_{n-1})\in \Phi^{j-1}_{n-1}\},\\
		\mathcal{B}^{j-1}_{n-2}&= \{(x_2,x_3,\dots,x_{n-2},0,2), (x_2,x_3,\dots,x_{n-2},1,1):(x_2,x_3,\dots,x_{n-2})\in \Phi^{j-1}_{n-2}\},\\
		\mathcal{B}^{j-1}_{n-3}&= \{(x_2,x_3,\dots,x_{n-3},0,0,2), (x_2,x_3,\dots,x_{n-3},0,1,1), \\
		&\hspace{3cm} (x_2,x_3,\dots,x_{n-3},1,0,1):(x_2,x_3,\dots,x_{n-3}) \in \Phi^{j-1}_{n-3}\},\\
		&\ \ \vdots\\
		\mathcal{B}^{j-1}_{j-1}&= \{(x_2,x_3,\dots,x_{j-1},0,0,\dots,0,2), (x_2,x_3,\dots,x_{j-1},0,0,\dots,0,1,1),\dots,\\
		&\ \ (x_2,x_3,\dots,x_{j-1},0,1,0\dots,0,1), (x_2,x_3,\dots,x_{j-1},1,0,\dots,0, 1):(x_2,x_3,\dots,x_{j-1})\in \Phi^{j-1}_{j-1}\}.
\end{align*}}

\begin{remark}
For $n\ge 1$, the cardinality of $\Phi^j_n$, $j\ge 1$ is a polynomial in $n$ of order $2j-1$. For different values of $j$, these are given by
{\small\begin{align*}
	|\Phi^1_n|&=n,\\
	|\Phi^2_n|&=	|\Phi^1_n|-1+\sum_{r_0=1}^{n-2}r_0|\Phi^1_{n-r_0}|=n-1+\sum_{r_0=1}^{n-2}r_0(n-r_0),\\
	|\Phi^3_n|&=\sum_{r_0=2}^{n-1}(n-r_0)\Big(r_0-1+\sum_{r_1=1}^{r_0-2}r_1(r_0-r_1)\Big),\\
	|\Phi^4_n|&=\sum_{r_0=3}^{n-1}\sum_{r_1=2}^{r_0-1}(n-r_0)(r_0-r_1)\Big(r_1-1+\sum_{r_2=1}^{r_1-2}r_2(r_1-r_2)\Big),\\
	&\ \  \vdots\\
	|\Phi^j_n|&=\sum_{r_0=j-1}^{n-1}\sum_{r_1=j-2}^{r_0-1}\sum_{r_2=j-3}^{r_1-1}\dots\sum_{r_{j-3}=2}^{r_{j-4}-1}(n-r_0)\Big(r_{j-3}-1+\sum_{r_{j-2}=1}^{r_{j-3}-2}r_{j-2}(r_{j-3}-r_{j-2})\Big)\prod_{i=0}^{j-4}(r_i-r_{i+1}).
\end{align*}	}
\end{remark}
\begin{remark}
Equivalently, we have the following explicit form for its mean increments:
\begin{equation}\label{equie_xn+1}
	\mathbb{E}(X_{n+1})=\sum_{j=1}^{n}\alpha^j\sum_{k=2(j-1)}^{2j}\sum_{\Omega_{k,n}^j}\prod_{i=2}^{n}\Big(\frac{1}{i}\Big)^{x_i}, \ n\ge1.
\end{equation}
Here, $\Omega_{k,n}^j=\{(x_2,x_3,\dots,x_n):x_i\in\{0,1,2\}, \ x_n\ne0,\ x_2+x_3+\dots+x_n=k,\ x_i\text{'s satisfy C1}\}$, where the condition C1 is as follows: \paragraph{C1} For each $2\le i\le n$ such that $x_i=2$, the cardinality of the set $D(i)=\{j: i< j\le n,\  x_j=1\}$ is even, that is,  $\#D(i)=2k$, $k\ge0$. 
\end{remark}

\begin{remark}
	Let us consider
	\begin{equation*}
		X_1=\begin{cases}
			+1 \ \  \text{with probability $q\in(0,1)$},\\
			-1 \ \ \text{with probability $1-q$}.
		\end{cases}
	\end{equation*}
	Then,
	$\mathbb{E}(X_1)=2q-1=\beta$ (say), $\mathbb{E}(X_2)=\beta\alpha$
	and from \eqref{Ex_n+1}, we have
	\begin{equation*}
		\mathbb{E}(X_{n+1})=\beta \sum_{j=1}^{n}\alpha^j\sum_{(x_2,x_3,\dots,x_n)\in\Phi^j_n}\frac{1}{2^{x_2}3^{x_3}\dots n^{x_n}},\ n\ge2.
	\end{equation*}
 Equivalently, from \eqref{equie_xn+1}, we get
	\begin{equation*}
	\mathbb{E}(X_{n+1})=\beta \sum_{j=1}^{n}\alpha^j\sum_{k=2(j-1)}^{2j}\sum_{\Omega_{k,n}^j}\prod_{i=2}^{n}\Big(\frac{1}{i}\Big)^{x_i}, \ n\ge1.
	\end{equation*}
\end{remark}
\begin{remark}
	For $n\ge0$, the mean increment $\mathbb{E}(X_{n+1})$ is a polynomial in $\alpha$ of degree $n$. 
\end{remark}

Next, we obtain the expressions for the mean displacement of the walker performing ERW with random memory. 
\subsection{Mean displacement}
From $\eqref{meanSn}$, we have
{\small\begin{align*}
	\mathbb{E}(S_1)&=\mathbb{E}(X_1)=1,\\
		\mathbb{E}(S_2)&=\mathbb{E}(S_1)+\alpha\mathbb{E}(S_1)=1+\alpha,\\
		\mathbb{E}(S_3)&=\mathbb{E}(S_2)+\frac{\alpha}{2}\Big(\mathbb{E}(S_1)+\frac{\mathbb{E}(S_2)}{2}\Big)=1+\alpha\Big(1+\frac{1}{2}+\frac{1}{2^2}\Big)+\frac{\alpha^2}{2^2}= \sum_{j=0}^{2}\alpha^j\sum_{x_2\in\Theta^j_2}\frac{1}{2^{x_2}},
	\end{align*}
	where $\Theta^0_2=\{0\}$, $\Theta^1_2=\{0, 1, 2\}$ and $\Theta^2_2=\{2\}$. 
		Also, 
	\begin{align*}
		\mathbb{E}(S_4)&=\mathbb{E}(S_3)+\frac{\alpha}{3}\Big(\mathbb{E}(S_1)+\frac{\mathbb{E}(S_2)}{2}+\frac{\mathbb{E}(S_3)}{3}\Big)\\
		&=1+\alpha\Big(1+\frac{1}{2}+\frac{1}{2^2}+\frac{1}{3}+\frac{1}{2\cdot3}+\frac{1}{3^2}\Big) +\alpha^2\Big(\frac{1}{2^2}+\frac{1}{2\cdot3}+\frac{1}{3^2}+\frac{1}{2\cdot3^2}+\frac{1}{2^2\cdot3^2}\Big) +\frac{\alpha^3}{2^2\cdot3^2}\\
		&=\sum_{j=0}^{3}\alpha^j\sum_{(x_2,x_3)\in\Theta^j_3}\frac{1}{2^{x_2}3^{x_3}},
	\end{align*}
	where $\Theta^0_3=\{(0,0)\}$,  $\Theta^1_3=\{(0,0),(1,0),(2,0),(0,1),(1,1),(0,2)\}$, $\Theta^2_3=\{(2,0),(1,1),(0,2),(1,2)$, $(2,2)\}$ and $\Theta^3_3=\{(2,2)\}$. 
	
	Similarly, we have
	\begin{align*}
		\mathbb{E}(S_5)&= \sum_{j=0}^{4}\alpha^j\sum_{(x_2,x_3,x_4)\in\Theta^j_4}\frac{1}{2^{x_2}3^{x_3}4^{x_4}},
	\end{align*}
	where
	\begin{align*}
		\Theta^0_4&=\{(0,0,0)\},\\
		\Theta^1_4&=\{(0,0,0),(1,0,0),(2,0,0),(0,1,0),(1,1,0),(0,2,0),(0,0,1),(1,0,1),(0,1,1),(0,0,2)\},\\
		\Theta^2_4&=\{(2,0,0),(1,1,0),(0,2,0),(1,2,0),(2,2,0),\\
		&\hspace{.7cm} (1,0,1),(0,1,1),(1,1,1),(2,1,1),(0,0,2),(1,0,2), (2,0,2),(0,1,2),(1,1,2),(0,2,2)\},\\
		\Theta^3_4&=\{(2,2,0),(2,1,1),(2,0,2),(1,1,2),(0,2,2),(1,2,2),(2,2,2)\},\\
		\Theta^4_4&=\{(2,2,2)\}.
	\end{align*}

	Proceeding inductively, we get the mean displacement of ERW with random memory in the following form:
\begin{equation}\label{ES_n}
\mathbb{E}(S_{n+1})=\sum_{j=0}^{n}\alpha^j\sum_{(x_2,x_3,\dots,x_n)\in\Theta^j_n}\frac{1}{2^{x_2}3^{x_3}\dots n^{x_n}}, \ n\ge2.
\end{equation}}
Here, $\Theta^0_n=\{(0,0,\dots,0)\}\subset\mathbb{R}^{n-1}$ and  $\Theta^j_n=\Lambda^{j}_{n-1}\cup\Lambda^{j-1}_{n-1}\cup\Lambda^{j-1}_{n-2}\cup\dots\cup\Lambda^{j-1}_{j-1}$, where
\begin{align*} \Lambda^{j}_{n-1}&=\{(x_2,x_3,\dots,x_{n-1},0):(x_2,x_3,\dots,x_{n-1})\in \Theta^j_{n-1}\},\\
	\Lambda^{j-1}_{n-1}&= \{(x_2,x_3,\dots,x_{n-1},2):(x_2,x_3,\dots,x_{n-1})\in \Theta^{j-1}_{n-1}\},\\
	\Lambda^{j-1}_{n-2}&= \{(x_2,x_3,\dots,x_{n-2},1,1):(x_2,x_3,\dots,x_{n-2})\in \Theta^{j-1}_{n-2}\},\\
	\Lambda^{j-1}_{n-3}&= \{(x_2,x_3,\dots,x_{n-3},1,0,1):(x_2,x_3,\dots,x_{n-3}) \in \Theta^{j-1}_{n-3}\},\\
	&\ \ \vdots\\
	\Lambda^{j-1}_{j-1}&= \{(x_2,x_3,\dots,x_{j-1},1,0,\dots,0, 1):(x_2,x_3,\dots,x_{j-1})\in \Theta^{j-1}_{j-1}\},
\end{align*}
such that
\begin{equation*}
\Lambda^{j-1}_{n-k}=\begin{cases}
\{(0,0,\dots,0,1)\}\subset\mathbb{R}^{k-1},\ n=k,\ j=1,\\
\{(1,0,\dots,0,1)\}\subset\mathbb{R}^{k},\ n=k+1,\ j\in\{1,2\}.
\end{cases}
\end{equation*}

\begin{remark}
For $n\ge 1$, the cardinality of $\theta^j_n$, $j\ge 1$ is a polynomial in $n$ of order $2j$. The cardinalities for different values of $j$ are given by
\begin{align*}
	|\theta^1_n|&=\sum_{r_0=1}^{n}r_0=\frac{n(n+1)}{2},\\
	|\theta^2_n|&=\sum_{r_0=1}^{n-1}(n-r_0)|\theta^1_{r_0}|,\\
	|\theta^3_n|&=\sum_{r_0=2}^{n-1}(n-r_0)\sum_{r_1=1}^{r_0-1}(r_0-r_1)|\theta^1_{r_1}|=\sum_{r_0=1}^{n-1}\sum_{r_1=1}^{r_0-1}(n-r_0)(r_0-r_1)|\theta^1_{r_1}|,\\
	|\theta^4_n|&=\sum_{r_0=3}^{n-1}\sum_{r_1=2}^{r_0-1}\sum_{r_2=1}^{r_1-1}(n-r_0)(r_0-r_1)(r_1-r_2)|\theta^1_{r_{2}}|,\\
	&\ \  \vdots\\
	|\theta^j_n|&=\sum_{r_0=j-1}^{n-1}\sum_{r_1=j-2}^{r_0-1}\sum_{r_2=j-3}^{r_1-1}\dots\sum_{r_{j-2}=1}^{r_{j-3}-1}(n-r_0)\Big(\prod_{i=0}^{j-3}(r_i-r_{i+1})\Big)|\theta^1_{r_{j-2}}|.
\end{align*}	
\end{remark}
\begin{remark}
Equivalently, the mean displacement of ERW with random memory has the following explicit form:
\begin{equation}\label{labES_n}
\mathbb{E}(S_{n+1})=1+\sum_{j=1}^{n}\alpha^j\sum_{k=2(j-1)}^{2j}\sum_{\Psi^j_{k,n}}\prod_{i=2}^{n}\Big(\frac{1}{i}\Big)^{x_i},\ n\ge 1.
\end{equation} 
Here, $\Psi^j_{k,n}=\{(x_2,x_3,\dots,x_n): x_i\in\{0,1,2\}, \ x_2+x_3+\dots+x_n=k, \ x_i\text{'s} \ \text{satisfy C1}\}$.
\end{remark}
\begin{remark}
	Let us consider
	\begin{equation*}
		X_1=\begin{cases}
			+1 \ \  \text{with probability $q\in(0,1)$},\\
			-1 \ \ \text{with probability $1-q$}.
		\end{cases}
	\end{equation*}
Then,
	$\mathbb{E}(S_1)=2q-1=\beta$ (say), $\mathbb{E}(S_2)=\beta\alpha$
	and from \eqref{ES_n}, we get
	\begin{equation*}
		\mathbb{E}(S_{n+1})=\beta\sum_{j=0}^{n}\alpha^j\sum_{(x_2,x_3,\dots,x_n)\in\Theta^j_n}\frac{1}{2^{x_2}3^{x_3}\dots n^{x_n}},\ n\ge2.
	\end{equation*}
Equivalently, from \eqref{labES_n}, we have
\begin{equation*}
\mathbb{E}(S_{n+1})=\beta+ \beta\sum_{j=1}^{n}\alpha^j\sum_{k=2(j-1)}^{2j}\sum_{\Psi^j_{k,n}}\prod_{i=2}^{n}\Big(\frac{1}{i}\Big)^{x_i},\ n\ge 1.
\end{equation*}
\end{remark}
\begin{remark}	
	For $n\ge0$, the mean displacement  $\mathbb{E}(S_{n+1})$ is a polynomial in $\alpha$ of degree $n$. 
\end{remark}

\section*{Acknowledgement}
The first author thanks Government of India for the grant of Prime Minister's Research Fellowship, ID 1003066.

\end{document}